\documentclass[12pt]{amsart}

\usepackage{mathptmx}      
\usepackage{helvet}            
\usepackage{courier}           
\usepackage{graphicx}              
\usepackage[bottom]{footmisc}
\usepackage{color}
\usepackage{mathrsfs}%
\usepackage{amsmath}%
\usepackage{amsfonts}%
\usepackage{amssymb}%
\usepackage{cite}
\newtheorem{theorem}{Theorem}[section]

\newtheorem{PLprop}[theorem]{Proposition}

\begin{document}
\title{The time singular limit for a fourth-order\\ damped wave equation for MEMS}
\author{Philippe Lauren\c{c}ot}
\address{Institut de Math\'ematiques de Toulouse, UMR~5219, Universit\'e de Toulouse, CNRS \\ F--31062 Toulouse Cedex 9, France}
\email{laurenco@math.univ-toulouse.fr}
\author{Christoph Walker}
\address{Leibniz Universit\"at Hannover, Institut f\" ur Angewandte Mathematik, Welfengarten 1, D--30167 Hannover, Germany}
\email{walker@ifam.uni-hannover.de}
%
%

\begin{abstract}
We consider a free boundary problem modeling electrostatic microelectromechanical systems. The model consists of a fourth-order damped wave equation for the elastic plate displacement which is coupled to an elliptic equation for the electrostatic potential. We first review some recent results on existence and non-existence of steady-states as well as on local and global well-posedness of the dynamical problem, the main focus being on the possible touchdown behavior of the elastic plate. We then investigate the behavior of the solutions in the time singular limit when the ratio between inertial and damping effects tends to zero.
\end{abstract}

\maketitle

\pagestyle{myheadings}
\markboth{\sc{Ph.  Lauren\c cot \& Ch. Walker}}{\sc{Time singular limit for the damped MEMS wave equation}}

\section{Introduction} \label{LWsec:int}

\begin{figure}
\centering\includegraphics[width=10cm]{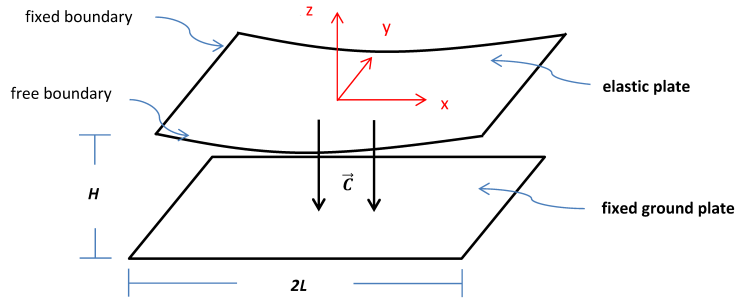}
\caption{\small Sketch of an idealized electrostatic MEMS device.}\label{PLMEMS1b}
\end{figure}

An idealized electostatically actuated microelectromechanical system (MEMS) consists of a fixed horizontal ground plate held at zero potential above which an elastic plate (or membrane) held at potential $V$ is suspended, see Figure~\ref{PLMEMS1b}.
A Coulomb force is generated by the potential difference across the device and results in a displacement of the elastic plate, thereby converting  electrostatic energy into mechanical energy, see \cite{PLEGG10,PLPB03} for a more detailed account and further references. After a suitable scaling and assuming homogeneity in the transversal horizontal direction (i.e. no $y$-dependence in Figure~\ref{PLMEMS1b}), the ground plate is assumed to be located at $z=-1$ and the plate displacement $u=u(t,x)\in (-1,\infty)$ evolves according to
\begin{equation}\label{PLueq}
\begin{split}
\gamma^2\partial_t^2u+\partial_tu & +\beta \partial_x^4 u - \tau \partial_x^2 u \\
 & = - \lambda \left( \varepsilon^2\ |\partial_x\psi(t,x,u(t,x))|^2 + |\partial_z\psi(t,x,u(t,x))|^2 \right)\ 
 \end{split}
\end{equation}
for $t>0$ and $x\in I:=(-1,1)$ with clamped boundary conditions
\begin{equation}
u(t,\pm 1) = \partial_x u(t,\pm 1) = 0\ , \quad t>0\ ,\label{PLubc}
\end{equation}
and initial conditions
\begin{equation}
u(0,x) = u^0(x)\ ,\quad \gamma^2\partial_t u(0,x) = \gamma^2 u^1(x)\ ,\quad x\in I \label{PLuic}\ .
\end{equation}
In \eqref{PLueq}, $\gamma^2\ge 0$ measures the ratio of inertial and damping forces which are given by the second and first order time derivatives, respectively, while $\beta\partial_x^4 u$ with \mbox{$\beta\ge 0$} and $-\tau\partial_x^2 u$ with $\tau\ge 0$ account for bending and stretching of the elastic plate, respectively. The right hand side of \eqref{PLueq} reflects the electrostatic forces exerted on the elastic plate, where the parameter $\lambda>0$ is proportional to the square of the voltage difference  between the two components, and the parameter $\varepsilon>0$ denotes the aspect ratio of the device (that is, the ratio height/length). The boundary conditions \eqref{PLubc} describe an elastic plate being clamped at its fixed boundary. Finally, the electrostatic potential $\psi=\psi(t,x,z)$ satisfies a rescaled Laplace equation in the  time-varying region 
$$
\Omega(u(t)) := \left\{ (x,z)\in I\times (-1,\infty)\ :\ -1 < z < u(t,x) \right\}
$$ 
between the ground plate and the elastic plate which reads
\begin{align}
\varepsilon^2 \,\partial_x^2\psi+\partial_z^2\psi & = 0 \ ,\quad &&(x,z)\in \Omega(u(t))\ ,\quad t>0\ , \label{PLpsieq} \\
\psi(t,x,z) & = \frac{1+z}{1+u(t,x)}\ ,\quad &&(x,z)\in\partial\Omega(u(t))\ ,\quad t>0\ .\label{PLpsibc}
\end{align}
Note that, when $\gamma>0$, equation~\eqref{PLueq} is a hyperbolic nonlocal semilinear fourth-order equation for the plate displacement $u$, which is coupled to the second-order  elliptic equation \eqref{PLpsieq} in the moving domain $\Omega(u(t))$ for the electrostatic potential $\psi$. If damping effects dominate over inertia effects one may set $\gamma=0$ in \eqref{PLueq}-\eqref{PLuic} and thus obtains a parabolic equation for $u$. 

A noteworthy feature of the above model is that it is only meaningful as long as the elastic plate does not touch down on the ground plate, that is, the deflection $u$ satisfies $u>-1$. From a physical point of view it is expected that above a certain critical threshold of $\lambda$, the elastic plate ``pulls in'' and smashes down on the ground plate. Obviously, the stable operating conditions of a given MEMS device heavily depend on the possible occurrence of this so-called {\em ``pull in'' instability}. 
Mathematically, the touchdown singularity manifests in the definition of $\Omega(u(t))$ which becomes disconnected if $u(t,x)$ reaches the value $-1$ at some point $x$, but also in the right hand side of \eqref{PLueq} as $\partial_z\psi$ becomes singular at such points since $\psi=1$ along $z=u$ while $\psi=0$ along $z=-1$.

\subsection{State of the Art}

According to the previous discussion, the mathematical investigation aims at showing that the parameter $\lambda$ indeed governs the dynamics of \eqref{PLueq}-\eqref{PLpsibc}, in particular the touchdown behavior and the closely related issues of global well-posedness and existence of steady states. More precisely, above a certain threshold value of $\lambda$ it is conjectured that solutions to \eqref{PLueq}-\eqref{PLpsibc}  cease to exist globally in time and that there are no steady-states, while for $\lambda$ below this critical value, solutions are global and there are at least two steady states. Moreover, if solutions do not exist globally, then the elastic plate pulls in at some finite time $T_c<\infty$, i.e., 
\begin{equation}
\lim_{t\to T_c} \min_{x\in I}\{ u(t,x) \} = -1\ . \label{PLtouchdown}
\end{equation}
In the special situation of the so-called {\em small aspect ratio model} which corresponds to setting \mbox{$\varepsilon=0$}  in \eqref{PLueq}-\eqref{PLpsibc}, it turns out that  the electrostatic potential $\psi$ is explicitly given by $\psi(t,x,z) = (1+z)/(1+u(t,x))$ in $\Omega(u(t))$ and the full system \eqref{PLueq}-\eqref{PLpsibc} reduces to a singular evolution equation only involving~$u$. In this case, a quite complete characterization of the expected dynamics -- confirming almost all of these conjectures -- is obtained in \cite{PLGLY14, PLLWxyy}, see also \cite{PLGu10, PLLL12} for further information as well as \cite{PLEGG10} and the references therein for the small aspect ratio model in general.

In contrast, due to the present coupling, the free boundary problem with $\varepsilon>0$ turns out to be even more involved and the literature is far more sparse in this case. A series of recent papers, however, addresses these questions for the free boundary problem when $\gamma=\beta=0$: see \cite{PLLaurencotWalker_ARMA} for steady-state solutions, \cite{PLELW1} for the corresponding parabolic problem, and \cite{PLELWyy,PLELW3} for a quasilinear version thereof. The fourth-order case $\beta>0$ with $\gamma\ge 0$ is investigated in \cite{PLLWxx} and we shall review its main results below (see also \cite{PLLWxx2} for a quasilinear version for $\beta>0$ and $\gamma=0$).

In all the just cited references on the free boundary problem, a very crucial ingredient in the analysis is the understanding of the elliptic problem \eqref{PLpsieq}-\eqref{PLpsibc} in the domain $\Omega(u)$ in dependence of a given (free) boundary described by a function $u:[-1,1]\rightarrow (-1,\infty)$ for a fixed time~$t$ (being suppressed for the moment). In particular, precise information on the gradient trace of the potential~$\psi=\psi_u$ on the elastic plate is required as a function of $u$. For this, one can transform the Laplace equation \eqref{PLpsieq}-\eqref{PLpsibc} for $\psi_u$ to an elliptic problem in the
fixed rectangle $I\times (0,1)$ (with coefficients depending on $u$ and its $x$-derivatives up to second order and being singular in case that $u$ approaches $-1$) for a transformed electrostatic potential $\phi_u$ given by
\begin{equation}\label{PLphi}
\phi_u(x,\eta) := \psi_u\big( x, (1+u(x))\eta - 1 \big) \ , \quad (x,\eta)\in I\times (0,1)\ .
\end{equation}
Using then elliptic regularity theory and pointwise multiplications in Sobolev spaces, the following key result can be shown \cite[Proposition~5]{PLELW1}:

\begin{PLprop} \label{PLP1a}
Given $2\alpha\in (0,1/2)$ and $\kappa\in (0,1)$ define an  open subset of $H^{2+2\alpha}(I)$ by
\begin{equation*}
\begin{split}
S_{\alpha}(\kappa):=\big\{v\in H^{2+2\alpha}(I)\,;\, & v(\pm 1)=\partial_xv(\pm 1)=0\,,\,\|v\|_{H^{2+2\alpha}(I)}< 1/\kappa, \\
&  \;\;\text{ and }\; v(x)>-1+\kappa\, \text{ for } x\in I \big\} \ . 
\end{split}
\end{equation*}
Then, for each $u\in S_{\alpha}(\kappa)$, there is a unique solution $\psi=\psi_u\in H^2(\Omega(u))$ to \eqref{PLpsieq}-\eqref{PLpsibc} and
\begin{equation}\label{PLest}
\|\psi_u\|_{H^2(\Omega(u))}\le C_0(\kappa)\ ,\quad u\in S_{\alpha}(\kappa)\ ,
\end{equation}
with 
\begin{equation}\label{PL10A}
\|\phi_{u_1}-\phi_{u_2}\|_{H^2(I\times (0,1))}\le C_0(\kappa) \|u_1-u_2\|_{H^{2+2\alpha}(I)}\ ,\quad u_1, u_2\in S_{\alpha}(\kappa)\ ,
\end{equation}
for some positive constant $C_0(\kappa)$ depending on $\kappa$ and also on $\alpha$, $\beta$, $\tau$, and $\varepsilon$, but not on $u\in S_{\alpha}(\kappa)$. Moreover, the mapping 
\begin{equation*}
g:S_\alpha(\kappa)\rightarrow H^{2\alpha}(I)\ ,\quad u\mapsto \varepsilon^2 |\partial_x \psi_u(x,u(x))|^2 + |\partial_z \psi_u(x,u(x))|^2 
\end{equation*}
is analytic, bounded, and uniformly Lipschitz continuous. 
\end{PLprop}


In fact, estimate~\eqref{PLest} follows from \cite[Lemma~6]{PLELW1} while \eqref{PL10A} is shown in  \cite[Eq.~(38)]{PLELW1}. Importantly, the minimal regularity required in order to control the potential $\psi_u$ in terms of suitable norms of $u$ appears to be that the latter belongs to $W_q^2(I)$ for some $q>2$, whence necessarily $2\alpha>0$ above (see \cite[Proposition~5]{PLELW1} for a more precise result). The regularity properties of the map $g$ are stated in this form for simplicity but are still sufficient in the fourth-order case $\beta>0$ considered in the following. 

As pointed out above, this case is investigated in \cite{PLLWxx}. In particular, existence and non-existence of steady states are derived in \cite{PLLWxx} in dependence of the voltage value~$\lambda$. While the former is a rather immediate consequence of the implicit function theorem (once Proposition~\ref{PLP1a} is established), the latter is based on a nonlinear version of the eigenfunction method which involves a positive eigenfunction in $H^4(I)$ associated to a positive eigenvalue of the fourth-order operator $\beta\partial_x^4-\tau\partial_x^2$ subject to clamped boundary conditions \cite{PLGr02,PLLWxy,PLOw97}. The result reads \cite[Theorem~1.7]{PLLWxx}:

\begin{theorem}[Steady States]\label{PLTStable2}
\begin{itemize}
\item[(i)] {\bf Existence:} There is $\lambda_s>0$ such that for each $\lambda\in (0,\lambda_s)$, there exists an asymptotically stable steady state $(U_\lambda,\Psi_\lambda)$ to \eqref{PLueq}-\eqref{PLpsibc} with $U_\lambda\in H^4(I)$ satisfying $-1<U_\lambda< 0$ in $I$ and  $\Psi_\lambda\in H^2(\Omega(U_\lambda))$.
\item[(ii)]\ {\bf Non-Existence:} There is ${\lambda_c}\ge \lambda_s$ such that there is no (sufficiently smooth) steady state $(u,\psi)$ to \eqref{PLueq}-\eqref{PLpsibc} for $\lambda>{\lambda_c}$.
\end{itemize}
\end{theorem}

Yet open problems are whether ${\lambda_c}= \lambda_s$ and whether there is a second (unstable) steady state for $\lambda<\lambda_s$ as in the small aspect ratio model, see \cite{PLGLY14, PLLWxyy}.

In \cite{PLLWxx} also the well-posedness of the dynamical problem is addressed. Due to Proposition~\ref{PLP1a} one may write \eqref{PLueq}-\eqref{PLpsibc} as a single semilinear Cauchy problem for the plate displacement $u$ (and its time derivative) that one can then solve by means of semigroup theory. We recall here the main statements from \cite[Propositions~3.1 \&~3.2, Corollaries~5.7 \&~5.10]{PLLWxx} and indicate an explicit dependence on the parameter $\gamma$ for future purposes:

\begin{theorem}[Well-Posedness]\label{PLT11}
Let $\gamma\ge 0$ and $2\alpha\in (0,1/2)$. Consider an initial condition $(u^0,u^1)$ in $H^{4+2\alpha}(I)\times H^{2+2\alpha}(I)$  satisfying $u^0(\pm 1) = \partial_x u^0(\pm 1)=u^1(\pm 1) = \partial_x u^1(\pm 1) = 0$ and such that $u^0> -1$ in~$I$. Then the following hold:

\begin{itemize}
\item[(i)] {\bf Local Existence:} For each $\lambda>0$, there is a unique solution $(u_\gamma,\psi_\gamma)$ to \eqref{PLueq}-\eqref{PLpsibc} on the maximal interval of existence $[0,T_{\gamma})$ in the sense that 
$$
\begin{array}{lcl}
u_\gamma\in C\big([0,T_\gamma),H^{2+2\alpha}(I)\big)\cap C^1\big([0,T_\gamma),H^{2\alpha}(I)\big) & \text{ if } & \gamma>0 \ , \\
 & & \\
u_0\in C\big([0,T_\gamma),H^{4}(I)\big)\cap C^1\big([0,T_\gamma),L_2(I)\big) & \text{ if } & \gamma=0 \ , 
\end{array}
$$
with
$$ 
\partial_t^k u_\gamma\in L_1\big((0,T), H^{4+2\alpha-2k}(I)\big)\ ,\quad k=0,1,2\ ,\quad T\in (0,T_\gamma)\ , \quad  \gamma>0\ ,
$$
satisfies \eqref{PLueq}-\eqref{PLuic} together with
$$
u_\gamma(t,x)>-1\ ,\quad (t,x)\in [0,T_\gamma)\times I\ , 
$$ 
while $\psi_{u_\gamma}(t) \in H^2\big( \Omega(u_\gamma(t)) \big)$  solves \eqref{PLpsieq}-\eqref{PLpsibc} in $\Omega(u_\gamma(t))$ for each $t\in [0,T_\gamma)$.

\item[(ii)]\ {\bf Touchdown:} There is $\gamma_1>0$ such that, if $\gamma\in [0,\gamma_1]$, then the solution $(u_\gamma,\psi_\gamma)$ to \eqref{PLueq}-\eqref{PLpsibc} obeys the following criterion for global existence: if for each $T>0$ there is $\kappa(T)\in (0,1)$ such that 
$$
u_\gamma(t)\ge -1+\kappa(T)\ \text{ in } I
$$ 
for $t\in [0,T_\gamma)\cap [0,T]$, then $T_\gamma=\infty$.

\item[(iii)]\ \, {\bf Global Existence:} Given $\kappa\in (0,1)$, there are numbers $\lambda_*=\lambda_*(\gamma,\kappa)>0$ and $N_*=N_*(\gamma,\kappa)>0$ such that $T_\gamma=\infty$ provided that 
$$ 
\|(u^0,u^1)\|_{H^{4+2\alpha}(I)\times H^{2+2\alpha}(I)}\le N_*\ ,\quad u^0\ge -1+\kappa \ \text{ in }\ I\ ,
$$
and $\lambda\le \lambda_*$. In this case,  $u_\gamma\in L_\infty\big((0,\infty),H^{2+2\alpha}(I)\big)$ with
$$
\inf_{(t,x)\in [0,\infty)\times I} u_\gamma(t,x)>-1\ .
$$
\end{itemize}
\end{theorem}

Actually,  in case of the damping dominated limit $\gamma=0$, less regularity on the initial data is required while more regularity on the solution $(u_0,\psi_0)$ may be obtained, see \cite[Propositions 3.1 \& 3.6]{PLLWxx} for details. The global existence result for small $\lambda$ values stated in part (iii) of the above theorem is based on the exponential decay of the associated semigroup which stems from the damping term. This fact will be exploited further in Section~\ref{PLsec:MT}. Note that for small  voltage values $\lambda$, touchdown is impossible, even in infinite time. An interesting, but still lacking salient feature of the physical model is a relation between $\lambda_c\ge \lambda_s$ from Theorem~\ref{PLTStable2} and (an optimally chosen) $\lambda_*$. 

The probably most important contribution to be brought forward by Theorem~\ref{PLT11} is the global existence criterion stated in part~(ii) which implies that touchdown is the only singularity preventing global existence. This is in clear contrast to the second-order case $\beta=0$ considered in \cite{PLELW1,PLELWyy,PLELW3}, where --~in principle~-- a finite existence time $T_\gamma$ may also be due to a blowup of some Sobolev norm of $u_\gamma(t,\cdot)$ as $t\to T_\gamma$. Roughly speaking, this physically most relevant feature is achieved by fully exploiting the additional information coming from the fourth-order term as well as the underlying gradient flow structure of \eqref{PLueq}-\eqref{PLpsibc}, the latter seeming to have been unnoticed  so far though being inherent in the model derivation. Indeed, introducing the total energy 
\begin{equation*}
\mathcal{E}(u_\gamma) := \mathcal{E}_m(u_\gamma) - \lambda \mathcal{E}_e(u_\gamma) 
\end{equation*}
involving the mechanical energy 
\begin{equation*}
\mathcal{E}_m(u_\gamma) := \frac{\beta}{2} \|\partial_x^2 u_\gamma \|_{L_2(I)}^2+\frac{\tau}{2} \|\partial_x u_\gamma \|_{L_2(I)}^2 
\end{equation*}
and the electrostatic energy
\begin{equation}
\mathcal{E}_e(u_\gamma) :=   \int_{\Omega(u_\gamma)} \left[ \varepsilon^2 |\partial_x \psi_{u_\gamma}(x,z)|^2 + |\partial_z \psi_{u_\gamma}(x,z)|^2 \right]\ \mathrm{d}(x,z)\ , \label{PLElecEnergy}
\end{equation}
the following energy equality holds \cite[Propositions~1.3 \&~1.6]{PLLWxx}:

\begin{PLprop}[Energy Equality]\label{PLP11}
Under the assumptions of Theorem~\ref{PLT11}~(i),
\begin{equation}
\mathcal{E}(u_\gamma(t)) + \frac{\gamma^2}{2}\| \partial_t u_\gamma(t)\|_{L_2(I)}^2+\int_0^t \|\partial_t u_\gamma(s)\|_{L_2(I)}^2\ \mathrm{d}s = \mathcal{E}(u^0)+\frac{\gamma^2}{2}\| u^1\|_{L_2(I)}^2 \label{PLE}
\end{equation}
for $t\in [0,T_\gamma)$.
\end{PLprop}

Note, however, that the energy $\mathcal{E}$ is the sum of terms with different signs and is thus {\em not} coercive. The main difficulty in the proof of Proposition~\ref{PLP11} is the computation of the derivative of $\mathcal{E}_e(u_\gamma)$ with respect to $u_\gamma$ since its dependence  on $u_\gamma$ is somehow implicit and involves the domain $\Omega(u_\gamma)$. Nevertheless, the derivative can be interpreted as the shape derivative of the Dirichlet integral of $\psi_\gamma=\psi_{u_\gamma}$, which can be computed and shown to be equal to the right hand side of \eqref{PLueq} --~except for the sign~-- by shape optimization arguments \cite{PLLWxx}. An additional difficulty stems from the fact that the time regularity of $u_\gamma$ as stated in part (i) of Theorem~\ref{PLT11} is not sufficient for a direct computation and one rather has to use an approximation argument. 

To prove then the significant criterion for global existence from part~(ii) of Theorem~\ref{PLT11}, one may proceed as follows: As long as $u_\gamma(t,\cdot)$ stays away from $-1$, one may control the electrostatic energy $\mathcal{E}_e(u_\gamma(t))$ by the mechanical energy $\mathcal{E}_m(u_\gamma(t))$ and then derives from the time decrease of $\mathcal{E}(u_\gamma(t))$ implied by Proposition~\ref{PLP11} first a bound on the $H^2(I)$- norm of $u_\gamma(t)$ and subsequently also on higher  Sobolev norms by a bootstrapping argument which yields global existence. 

\subsection{The Time Singular Limit}

In many research papers -- mostly dedicated to the small aspect ratio model with $\varepsilon=0$ -- inertial effects are neglected from the outset as damping effects may be predominant, a few exceptions being \cite{PLGu10, PLKLNT11}. In this note we now shall investigate the behavior of the solutions in the damping dominated limit $\gamma^2\to 0$. Obviously, considering such a time singular limit from a mathematical point of view requires in particular a common interval of existence, independent of $\gamma$, that is, a lower bound on the maximal existence time~$T_\gamma$. This is provided by the first result of this paper:

\begin{PLprop}[Minimal Existence Time]\label{PLP2}
Let $2\alpha\in (0,1/2)$, $\gamma \in (0,\gamma_1]$, $\lambda>0$. Consider an initial condition $(u^0,u^1)\in H^{4+2\alpha}(I)\times H^{2+2\alpha}(I)$ such that $u^0\in S_\alpha(\kappa)$ for some $\kappa\in (0,1)$ and $u^1(\pm1)=\partial_x u^1(\pm 1)=0$. Let $(u_\gamma,\psi_\gamma)$ be the unique solution to \eqref{PLueq}-\eqref{PLpsibc} defined on the maximal interval of existence $[0,T_{\gamma})$. There are $\hat{\gamma} := \hat{\gamma}\big( \kappa,\gamma_1,\|u^1\|_{H_D^{2\alpha}(I)} \big) \in (0,\gamma_1]$, $N:=N(\kappa,\gamma_1)>0$, and $\Lambda := \Lambda(\kappa,\gamma_1)>0$ such that:
\begin{itemize}
\item[(i)] There is $\hat{T} :=\hat{T} \big(\lambda,\kappa,\gamma_1,\|u^0\|_{H^{4+2\alpha}(I)} \big) \in (0,\infty)$ such that, for all $\gamma\in (0, \hat{\gamma})$, \mbox{$T_\gamma > \hat{T}$} and $u_\gamma(t)\in S_\alpha(\kappa/2)$ for $t\in [0,\hat{T}]$.

\item[(ii)]\ If $\lambda\in (0,\Lambda)$ and $\|u^0\|_ {H^{4+2\alpha}(I)}\le N$, then $T_\gamma =\infty$ and $u_\gamma(t)\in S_\alpha(\kappa/2)$ for $t\ge 0$ and $\gamma\in (0,\hat{\gamma})$.
\end{itemize}
\end{PLprop}

The proof of this proposition is given in Section~\ref{PLsec:MT}. It relies on an exponential decay of the energy associated to the damped wave equation being independent of $\gamma\in [0, \gamma_1]$.

As a consequence we are in a position to investigate the damping dominated limit and prove that $(u_\gamma,\psi_\gamma)$ converges toward $(u_0,\psi_0)$ in a suitable sense as $\gamma^2\to 0$.

\begin{theorem}[Damping Dominated Limit]\label{PLT}
Under the assumptions of Proposition~\ref{PLP2}~(i) and as $\gamma^2\longrightarrow 0$,
\begin{equation}\label{PLs1}
u_\gamma\longrightarrow u_0\ \text{ in } \ C\big([0,\hat{T}],H^{2+2\xi}(I)\big)
\end{equation}
for each $\xi\in (0,\alpha)$
and
\begin{equation}\label{PLs2}
\phi_{u_\gamma}\longrightarrow \phi_{u_0} \ \text{ in } \ C\big([0,\hat{T}],H^{2}(I\times (0,1))\big)\ ,
\end{equation}
where $\phi_{u_\gamma}$ is the transformed electrostatic potential given by \eqref{PLphi} (with $u$ replaced by $u_\gamma$). In addition, if $u^0=u^1=0$, then 
\begin{equation}\label{PLs3}
\partial_t u_\gamma\longrightarrow \partial_t u_0\ \text{ in } \ L_p\big(0,\hat{T};H^{2\alpha}(I)\big)
\end{equation}
for each $p\in (1,\infty)$. Under the assumptions of Proposition~\ref{PLP2}~(ii), statements \eqref{PLs1}-\eqref{PLs3} are true for each $T>0$ instead of $\hat{T}$.
\end{theorem}

The proof of Theorem~\ref{PLT} is performed in Section~\ref{PLsec2}. It is based on compactness properties of $(u_\gamma,\psi_\gamma)_{\gamma\in (0,\gamma_1)}$ being provided by the energy functional $\mathcal{E}$.

\section{A Lower Bound on the Maximal Existence Time} \label{PLsec:MT}

In order to prove Proposition~\ref{PLP2}, we consider an initial condition $(u^0,u^1)$ belonging to $H_D^{4+2\alpha}(I)\times H_D^{2+2\alpha}(I)$ and such that $u^0\in S_\alpha(\kappa)$ for some $\kappa\in (0,1)$ and \mbox{$2\alpha\in (0,1/2)$}, where
$$
H_D^{\theta}(I):=\left\{\begin{array}{lll}
& \big\{v\in H^{\theta}(I)\,;\, v(\pm 1)=\partial_x v(\pm 1)=0\big\}\ , & \theta>\dfrac{3}{2}\ ,\\
& \big\{v\in H^{\theta}(I)\,;\, v(\pm 1)=0\big\}\ , & \dfrac{1}{2}<\theta<\dfrac{3}{2}\ ,\\
&   H^{\theta}(I)\ , & \theta<\dfrac{1}{2}\ .
\end{array}
\right.
$$
We fix $\gamma\in (0,\gamma_1]$ with $\gamma_1>0$ introduced in Theorem~\ref{PLT11}~(ii) and let $(u_\gamma,\psi_\gamma)$ with
$$
u_\gamma\in C\big([0,T_\gamma),H_D^{2+2\alpha}(I)\big)\cap C^1\big([0,T_\gamma),H_D^{2\alpha}(I)\big)$$
and
$$ 
\partial_t^k u_\gamma\in L_1\big((0,T), H_D^{4+2\alpha-2k}(I)\big)\ ,\quad k=0,1,2\ ,\quad T\in (0,T_\gamma)\ ,
$$
be the unique solution  to \eqref{PLueq}-\eqref{PLpsibc} on the maximal interval of existence $[0,T_{\gamma})$ as provided by Theorem~\ref{PLT11}. Then, introducing the operator
$$
A_\alpha :=\beta \partial_x^4 - \tau \partial_x^2\in \mathcal{L}\big( H_D^{4+2\alpha}(I),H_D^{2\alpha}(I)\big)
$$ 
we have
\begin{equation*}\label{PLuequ}
\gamma^2\frac{\mathrm{d}^2}{\mathrm{d} t^2} u_\gamma+\frac{\mathrm{d}}{\mathrm{d} t}u_\gamma+ A_\alpha u_\gamma = - \lambda g(u_\gamma)\ , \quad t\in (0,T_\gamma)\ ,
\end{equation*}
in $H_D^{2\alpha}(I)$, the function $g$ being defined in Proposition~\ref{PLP1a}.

We want to control a suitable norm of $u_\gamma(t)$ for which we basically use an idea from \cite[Section 2]{PLHaZa}, the difference mainly being the focus on estimates which are uniform with respect to $\gamma\in (0,\gamma_1]$. To this end, define 
$$
v(t):=u_\gamma(t)-u^0\quad \text{ and }\quad f(t):= -\lambda g(u_\gamma(t)) -A_\alpha u^0
$$ 
for $t\in [0,T_\gamma)$. Then $v$ solves the equation
\begin{equation}\label{PLm0}
\gamma^2\frac{\mathrm{d}^2}{\mathrm{d} t^2}v+\frac{\mathrm{d}}{\mathrm{d} t}v+ A_\alpha v  = f \ , \quad t\in (0,T_\gamma)\ ,
\end{equation}
in $H_D^{2\alpha}(I)$ with initial condition $(v(0),\partial_t v(0))=(0,u^1)$.
Recall that there are real numbers $c_2\ge c_1\ge c_0\ge 1$ such that
\begin{equation}\label{PLm1}
\|z\|_{H_D^{2\alpha}(I)}^2\le c_0 \|z\|_{H_D^{2+2\alpha}(I)}^2\le c_1 \|A_\alpha^{1/2}z\|_{H_D^{2\alpha}(I)}^2\le c_2 \|z\|_{H_D^{2+2\alpha}(I)}^2
\end{equation}
for all $z\in H_D^{2+2\alpha}(I)$. Then, defining
$$
E(t):= \left\|A_\alpha^{1/2} v(t)\right\|_{H_D^{2\alpha}(I)}^2+\gamma^2  \left\|\partial_t v(t)\right\|_{H_D^{2\alpha}(I)}^2\ ,\quad t\in (0,T_\gamma)\ ,
$$
and
$$
F(t):=\gamma\left\langle v(t)\,,\,\partial_t v(t) \right\rangle_{H_D^{2\alpha}(I)}\ ,\quad t\in (0,T_\gamma)\ ,
$$
we deduce from \eqref{PLm0}, \eqref{PLm1}, and the self-adjointness of $A_\alpha^{1/2}$ in $H_D^{2\alpha}(I)$ that 
\begin{equation}\label{PLm2a}
\frac{\mathrm{d}}{\mathrm{d} t} E(t)=-2\left\| \partial_t  v(t)\right\|_{H_D^{2\alpha}(I)}^2  +  2\left\langle f(t)\,,\,\partial_t v(t)\right\rangle_{H_D^{2\alpha}(I)}
\end{equation}
and 
\begin{equation}\label{PLm2b}
|F(t)| \le \frac{c_1}{2} E(t)
\end{equation}
for a. e. $t\in (0,T_\gamma)$. Next, let 
$$
b:=\min\left\{\frac{2}{2\gamma_1^2+c_1+1}\,,\,\frac{1}{2 c_1}\,,\,\frac{1}{\gamma_1 c_1}\right\}
$$ 
and introduce $G(t):=E(t)+b\gamma F(t)$ for $t\in (0,T_\gamma)$. According to \eqref{PLm1}-\eqref{PLm2b} and Young's inequality,
\begin{equation*}
\begin{split}
\frac{\mathrm{d}}{\mathrm{d} t} G(t)&= \left(-2+b\gamma^2\right) \left\| \partial_t v(t)\right\|_{H_D^{2\alpha}(I)}^2 -b\left\langle v(t)\,,\, \partial_t  v(t)+A_\alpha v(t)-f(t)\right\rangle_{H_D^{2\alpha}(I)}\\
&\quad +  2\left\langle f(t)\,,\, \partial_t  v(t) \right\rangle_{H_D^{2\alpha}(I)}\\ 
&\le \left(-2+b\gamma^2\right) \left\| \partial_t v(t)\right\|_{H_D^{2\alpha}(I)}^2 -b \left\|A_\alpha^{1/2}v(t)\right\|_{H_D^{2\alpha}(I)}^2 \\
&\quad +b\left(\frac{1}{4c_1}\|v(t)\|_{H_D^{2\alpha}(I)}^2+c_1 \left\| \partial_t v(t)\right\|_{H_D^{2\alpha}(I)}^2\right)\\
&\quad +\frac{b^2}{2} \|v(t)\|_{H_D^{2\alpha}(I)}^2 +\frac{1}{2} \|f(t)\|_{H_D^{2\alpha}(I)}^2+b \left\| \partial_t v(t)\right\|_{H_D^{2\alpha}(I)}^2 +\frac{1}{b}\|f(t)\|_{H_D^{2\alpha}(I)}^2\\
&\le\left(-2+b\gamma^2 +b c_1+b\right) \left\|\partial_t v(t)\right\|_{H_D^{2\alpha}(I)}^2 -\frac{b}{2}\left\|A_\alpha^{1/2}v(t)\right\|_{H_D^{2\alpha}(I)}^2\\
&\quad+\frac{b}{2}\left( \left(\frac{1}{2c_1} +b\right)\left\|v(t)\right\|_{H_D^{2\alpha}(I)}^2 -\left\|A_\alpha^{1/2}v(t)\right\|_{H_D^{2\alpha}(I)}^2\right)
\\
&\qquad +\left(\frac{1}{2}+\frac{1}{b}\right)\|f(t)\|_{H_D^{2\alpha}(I)}^2
\end{split}
\end{equation*}
for a. e. $t\in (0,T_\gamma)$. Since the choice of $b$ ensures that
$$
\frac{1}{2c_1} +b\le \frac{1}{c_1} \;\;\text{ and }\;\; - 2 + b (\gamma^2 + c_1 + 1) \le - b \gamma^2\ ,
$$
the third term in the right hand side is non-positive by \eqref{PLm1} and we obtain
\begin{equation*}
\frac{\mathrm{d}}{\mathrm{d} t} G(t)\le -\frac{b }{2} E(t) + \frac{b+2}{2b} \|f(t)\|_{H_D^{2\alpha}(I)}^2 \quad \text{ for a.e.\ } t\in (0,T_\gamma)\ .
\end{equation*}
Observe that \eqref{PLm2b} and the choice of $b$ also ensure
$$
\frac{1}{2}E(t)\le \left(1-\frac{c_1b\gamma}{2}\right) E(t)\le G(t)\le \left(1+\frac{c_1b\gamma}{2}\right) E(t) \le \left( 1 + \frac{c_1 b \gamma_1}{2} \right) E(t) 
$$
for $t\in (0,T_\gamma)$,
whence
$$
\frac{\mathrm{d}}{\mathrm{d} t} G(t)\le -\frac{b}{2+c_1 b \gamma_1} G(t) +  \frac{b+2}{2b}  \|f(t)\|_{H_D^{2\alpha}(I)}^2\quad  \text{for a.e.\ } t\in (0,T_\gamma)\ .
$$
Consequently, setting $\omega:=b/(2+c_1 b\gamma_1)$, 
\begin{equation*}
\begin{split} 
E(t)\le 2 G(t)\le \frac{b}{\omega} e^{-\omega t} E(0)  +\frac{b+2}{b \omega} \left(1-e^{-\omega t}\right) \sup_{s\in (0,t)}\left\{ \|f(s)\|_{H_D^{2\alpha}(I)}^2 \right\}
\end{split}
\end{equation*}
for $t\in (0,T_\gamma)$. Now, owing to \eqref{PLm1} and the definitions of $E(t)$ and $f(t)$, there is a constant $M:=M(\gamma_1)>0$ such that
\begin{equation*}\label{PLee}
\begin{split} 
\|&u_\gamma(t)-u^0\|_{H_D^{2+2\alpha}(I)}^2\\
&\ \le M \gamma^2 \|u^1\|_{H_D^{2\alpha}(I)}^2
+ M\left(1-e^{-\omega t}\right) \left[ \lambda^2\sup_{s\in (0,t)}\left\{ \|g(u_\gamma(s))\|_{H_D^{2\alpha}(I)}^2 \right\}+ \|u^0\|_{H_D^{4+2\alpha}(I)}^2\right]
\end{split}
\end{equation*}
for $t\in (0,T_\gamma)$. Since $u^0$ belongs to $S_{\alpha}(\kappa)$, it follows from its time continuity in $H_D^{2+2\alpha}(I)$ and the continuous embedding of $H_D^{2+2\alpha}(I)$ in $L_\infty(I)$ that
$$
\hat{T}_\gamma := \sup\left\{t_0\in (0,T_\gamma)\,;\, u_\gamma(t)\in S_\alpha(\kappa/2) \;\text{ for all }\; t\in [0,t_0)\right\} > 0\ .
$$
Then, by Proposition~\ref{PLP1a}, 
\begin{equation*}\label{PL3}
\| g(u_\gamma(t))\|_{H_D^{2\alpha}(I)}\le c_3(\kappa) := \sup_{w\in S_\alpha(\kappa/2)}\left\{ \| g(w)\|_{H_D^{2\alpha}(I)} \right\} \ ,\quad  t\in [0,\hat{T}_\gamma)\ ,
\end{equation*}
and we conclude that 
\begin{align*}
\|u_\gamma(t)-u^0\|_{H_D^{2+2\alpha}(I)}^2 & \le  M \gamma^2 \|u^1\|_{H_D^{2\alpha}(I)}^2 \\
& + M \left[ \lambda^2 c_3(\kappa)^2 + \|u^0\|_{H_D^{4+2\alpha}(I)}^2 \right] \left(1-e^{-\omega t}\right)
\end{align*}
 for $t\in [0,\hat{T}_\gamma)$. Therefore, since $u^0\in S_\alpha(\kappa)$ and since $H_D^{2+2\alpha}(I)$ embeds continuously in $L_\infty(I)$ with constant, say, $c_4\ge 1$, the previous inequality ensures that
\begin{align*}
\|u_\gamma(t)\|_{H_D^{2+2\alpha}(I)} & \le \|u_\gamma(t)-u^0\|_{H_D^{2+2\alpha}(I)} + \|u^0\|_{H_D^{2+2\alpha}(I)} < \frac{2}{\kappa}
\end{align*}
as soon as 
$$
M \gamma^2 \|u^1\|_{H_D^{2\alpha}(I)}^2 + M \left[ \lambda^2 c_3(\kappa)^2 + \|u^0\|_{H_D^{4+2\alpha}(I)}^2 \right] \left(1-e^{-\omega t}\right) < \frac{1}{\kappa^2}
$$
and 
\begin{align*}
u_\gamma(t) & \ge u^0 - \|u_\gamma(t)-u^0\|_{L_\infty(I)} \ge \kappa - 1 - c_4 \|u_\gamma(t)-u^0\|_{H_D^{2+2\alpha}(I)} > \frac{\kappa}{2} - 1
\end{align*}
as soon as 
$$
M \gamma^2 \|u^1\|_{H_D^{2\alpha}(I)}^2 + M \left[ \lambda^2 c_3(\kappa)^2 + \|u^0\|_{H_D^{4+2\alpha}(I)}^2 \right] \left(1-e^{-\omega t}\right) < \frac{\kappa^2}{4 c_4^2}\ .
$$
Thus, since $\kappa^2/(4 c_4^2) \le 1 \le 1/\kappa^2$, we deduce from the above analysis that $u_\gamma(t)$ belongs to $S_\alpha(\kappa/2)$ provided $t\in [0,\hat{T}_\gamma)$ and $\gamma\in (0,\gamma_1]$ satisfy
\begin{equation}
\gamma^2 \|u^1\|_{H_D^{2\alpha}(I)}^2 < \frac{\kappa^2}{8 M c_4^2} \label{PLvolvic}
\end{equation}
and
\begin{equation}
\left[ \lambda^2 c_3(\kappa)^2 + \|u^0\|_{H_D^{4+2\alpha}(I)}^2 \right] \left(1-e^{-\omega t}\right) < \frac{\kappa^2}{8 M c_4^2}\ . \label{PLevian}
\end{equation}
Therefore, there are 
$$
\hat{\gamma} := \hat{\gamma}\big( \kappa, \gamma_1, \|u^1\|_{H_D^{2\alpha}(I)} \big) \in (0,\gamma_1) \quad\text{ and }\quad \hat{T} := \hat{T}\big( \lambda, \kappa, \gamma, \|u^0\|_{H_D^{4+2\alpha}(I)} \big)> 0
$$ 
such that
$$
u_\gamma(t)\in S_\alpha(\kappa/2) \;\;\text{ for }\;\; t\in [0,\hat{T}_\gamma) \cap [0,\hat{T}] \;\;\text{ and }\;\; \gamma\in (0,\hat{\gamma})
$$
Recalling the definition of $\hat{T}_\gamma$, the previous statement implies in particular that \mbox{$\hat{T}_\gamma\ge \hat{T}$}. Finally, owing to the positivity of $\omega$, it is clear that if one requires that
$$
\lambda^2 c_3(\kappa)^2 + \|u^0\|_{H_D^{4+2\alpha}(I)}^2 < \frac{\kappa^2}{8 M c_4^2}\ 
$$
instead of \eqref{PLevian}, there are $\Lambda:=\Lambda(\kappa,\gamma_1)>0$ and $N:=N(\kappa,\gamma_1)>0$ such that $u_\gamma(t)$ belongs to $S_\alpha(\kappa/2)$ for all $t\in [0,T_\gamma)$ and $\gamma\in (0,\hat{\gamma})$ provided that $\lambda\in (0,\Lambda)$ and $\|u^0\|_{H_D^{4+2\alpha}(I)}\le N$, whence $T_\gamma=\infty$ by Theorem~\ref{PLT11}~(ii). This proves Proposition~\ref{PLP2}.


\section{The Time Singular Limit $\gamma^2\longrightarrow 0$}\label{PLsec2}

In order to prove Theorem~\ref{PLT} we stick to the notation from the previous section. Recall that 
\begin{equation}\label{PLm6}
u_\gamma(t)\in S_\alpha(\kappa/2)\ ,\quad t\in [0,\hat{T}]\ ,\quad \gamma\in (0,\hat{\gamma})\ .
\end{equation}
It then follows from \eqref{PLest} that
\begin{equation*}
\|\psi_{u_\gamma}(t)\|_{H^{2}(\Omega(u_\gamma(t)))}\le C_0(\kappa)\ ,\quad  t\in [0,\hat{T}]\ ,\quad \gamma\in (0,\hat{\gamma})\ .
\end{equation*}
This gives a uniform bound on the electrostatic energy $\mathcal{E}_e(u_\gamma(t))$ defined in \eqref{PLElecEnergy} so that \eqref{PLE} implies 
\begin{equation}\label{PLm8}
\frac{\gamma^2}{2}\|u_\gamma(t)\|_{L_2(I)}^2+\int_0^t\left\| \partial_t u_\gamma(s)\right\|_{L_2(I)}^2\,\mathrm{d} s\le c(\kappa)\ ,\quad  t\in [0,\hat{T}]\ ,\quad \gamma\in (0,\hat{\gamma})\ .
\end{equation}
Now, let $\xi\in (0,\alpha)$. Owing to \eqref{PLm6} and \eqref{PLm8}, the set $\{u_\gamma\,;\, \gamma\in (0,\hat{\gamma})\}$ is bounded in $L_\infty(0,\hat{T};H^{2+2\xi}(I))$ with $\{\partial_t u_\gamma(t)\,;\, \gamma\in (0,\hat{\gamma})\}$ bounded in $L_2\big( (0,\hat{T})\times I \big)$. We then infer from the compactness of the embedding of $H^{2+2\xi}(I)$ in $H^{2+2\alpha}(I)$ and \cite[Corollary~4]{PLSi87} that there are subsequence of $\gamma^2\longrightarrow 0$ (not relabeled) and $\bar u_0$ in $C\big([0,\hat{T}],H_D^{2+2\xi}(I)\big)$ such that
\begin{equation}\label{PLm9}
u_\gamma\longrightarrow \bar u_0\ \text{ in }\ C\big([0,\hat{T}],H_D^{2+2\xi}(I)\big)\ .
\end{equation}
Clearly, $\bar u_0(t)\in S_\alpha(\kappa/4)$ for $t\in [0,\hat{T}]$ by \eqref{PLm6} and \eqref{PLm9}. The latter and \eqref{PL10A} also imply
$$
\|\phi_{u_\gamma(t)}-\phi_{\bar u_0(t)}\|_{H^2(I\times (0,1))}\le C_0(\kappa/4) \|u_\gamma(t)-\bar u_0(t)\|_{H^{2+2\xi}(I)}
$$
for $t\in [0,\hat{T}]$ and $\gamma\in (0,\hat{\gamma})$. Consequently, Theorem~\ref{PLT} follows if we can show that $\bar u_0$ and $u_0$ coincide. To this end recall that the function $g:S(\kappa/4)\rightarrow H_D^{2\alpha}(I)$, defined in Proposition~\ref{PLP1a}, is uniformly Lipschitz continuous. In particular, from \eqref{PLm9} we deduce that for each $p\in (1,\infty)$, 
\begin{equation}\label{PL10}
g(u_\gamma)\longrightarrow g_0:=g(\bar u_0)\ \text{ in }\ L_p\big(0,\hat{T};H_D^{2\alpha}(I)\big)\ .
\end{equation}
Thus, if $v_\gamma$ denotes the solution to the linear Cauchy problem 
\begin{eqnarray*}
\gamma^2\frac{\mathrm{d}^2}{\mathrm{d} t^2}v+\frac{\mathrm{d}}{\mathrm{d} t}v+A_\alpha v & = &-\lambda g(u_\gamma)\ ,\quad t\in [0,\hat{T}]\ ,
\end{eqnarray*}
subject to zero initial conditions
$$
v(0)=\gamma^2 \partial_t v(0)=0 \ ,
$$
for $\gamma\in (0,\hat{\gamma})$ (with $v_0$ denoting accordingly the solution with $\gamma=0$), it follows from \eqref{PL10}, the fact that $-A_\alpha$ generates a strongly continuous cosine family in $H_D^{2\alpha}(I)$ as pointed out in \cite[Section 3.2]{PLLWxx}, and \cite[VI.Theorem~7.6]{PLFattorini}
that
\begin{equation}\label{PL110}
v_\gamma\longrightarrow v_0\ \text{ in }\ C\big([0,\hat{T}],H_D^{2\alpha}(I)\big)\ ,\quad \partial_t v_\gamma \longrightarrow\partial_t v_0 \ \text{ in }\  L_p\big(0,\hat{T};H_D^{2\alpha}(I)\big)\ .
\end{equation}
On the other hand, if $w_\gamma$ denotes the solution to the homogeneous Cauchy problem 
\begin{eqnarray*}
\gamma^2\frac{\mathrm{d}^2}{\mathrm{d} t^2}w+\frac{\mathrm{d}}{\mathrm{d} t}w+A_\alpha w & = &0\ ,\quad t>0\ ,
\end{eqnarray*}
subject to the initial conditions
$$
w(0)=u^0\ ,\quad \gamma^2 \partial_t w(0)=\gamma^2 u^1 \ ,
$$
for $\gamma\in (0,\hat{\gamma})$ (with $w_0$ denoting accordingly the solution with $\gamma=0$), then
\begin{equation}\label{PL1100}
w_\gamma\longrightarrow w_0\ \text{ in }\ C\big([0,\hat{T}],H_D^{2\alpha}(I)\big)
\end{equation}
owing to \cite[Theorem 3.2]{PLKisynski}. Clearly, by uniqueness of solutions to linear wave equations, we have $u_\gamma=v_\gamma+w_\gamma$, and consequently,  from \eqref{PLm9}, \eqref{PL110}, and \eqref{PL1100} we derive that $\bar u_0=v_0+w_0$ solves
$$
\frac{\mathrm{d}}{\mathrm{d} t} \bar u_0+A_\alpha \bar u_0=-\lambda g(\bar u_0)\ ,\quad t\in (0,\hat{T}]\ ,\qquad \bar u_0(0)=u^0\ .
$$
Since the above Cauchy problem has a unique solution according to Theorem~\ref{PLT11}, namely $u_0$ (restricted to $[0,\hat{T}]$), we conclude that $\bar u_0=u_0$
and since this limit is independent of the subsequence $\gamma^2\longrightarrow 0$, Theorem~\ref{PLT} is proven.


\end{document}